\documentclass[aip,graphicx]{revtex4-1}
\usepackage{amssymb}
\usepackage{amsmath}
\usepackage{graphicx}% Include figure files
\usepackage{psfrag}
\usepackage{lipsum}
\usepackage{caption}
\usepackage{xcolor}

\usepackage{dcolumn}% Align table columns on decimal point
\usepackage{bm}% bold math
\usepackage[utf8]{inputenc}
\usepackage[T1]{fontenc}
\usepackage{mathptmx} 

\usepackage{float}
\usepackage{tcolorbox}
\usepackage{hyperref}
\usepackage{subcaption}
\draft % marks overfull lines with a black rule on the right

\begin{document}

\title{Analysis of pseudo-spectral methods used for numerical simulation of turbulence}

\author{Tapan~K.~Sengupta}
\email[]{tksengupta@iitism.ac.in}
\affiliation{Dept. of Mechanical Engineering, Indian Institute of Technology (ISM) Dhanbad, Dhanbad-826004, India}

\author{Vajjala~K.~Suman}
%\email[]{vksuman@iitk.ac.in}
\affiliation{High Performance Computing Laboratory, Indian Institute of Technology Kanpur, Kanpur-208 016, India}
\affiliation{Computational \& Theoretical Fluid Dynamics Division, CSIR-NAL, Bangalore-560017, India}

\author{Prasannabalaji Sundaram}
%\email[]{prasanna@iitk.ac.in}
\affiliation{High Performance Computing Laboratory, Indian Institute of Technology Kanpur, Kanpur-208 016, India}

\author{Aditi Sengupta}
%\email[]{aditi@iitism.ac.in}
\affiliation{Dept. of Mechanical Engineering, Indian Institute of Technology (ISM) Dhanbad, Dhanbad-826004, India}

\date{\today}

\begin{abstract}
Global spectral analysis (GSA) is used as a tool to test the accuracy of numerical methods with the help of canonical problems of convection and convection-diffusion equation which admit exact solutions. Similarly, events in turbulent flows computed by direct numerical simulation (DNS) are often calibrated with theoretical results of homogeneous isotropic turbulence due to Kolmogorov, as given in {\it Turbulence} --U. Frisch, Cambridge Univ. Press, UK (1995). However, numerical methods for the simulation of this problem are not calibrated, as by using GSA of convection and/or convection-diffusion equation. This is with the exception in "A critical assessment of simulations for transitional and turbulence flows- Sengupta, T.K., In Proc. of IUTAM Symp. on {\it Advances in Computation, Modeling and Control of Transitional and Turbulent Flows}, pp 491-532, World Sci. Publ. Co. Pte. Ltd., Singapore (2016)", where such a calibration has been advocated with the help of convection equation. For turbulent flows, an extreme event is characterized by the presence of length scales smaller than the Kolmogorov length scale, a heuristic limit for the largest wavenumber present without being converted to heat. With growing computer power, recently many simulations have been reported using a pseudo-spectral method, with spatial discretization performed in Fourier spectral space and a two-stage, Runge-Kutta (RK2) method for time discretization. But no analyses are reported to ensure high accuracy of such simulations. Here, an analysis is reported for few multi-stage Runge-Kutta methods in the Fourier spectral framework for convection and convection-diffusion equations. We identify the major source of error for the RK2-Fourier spectral method using GSA and also show how to avoid this error and specify numerical parameters for achieving highest accuracy possible to capture extreme events in turbulent flows.

\end{abstract}


\maketitle %\maketitle must follow title, authors, abstract and \pacs

%\linenumbers

\section{Introduction}

%Introduce DNS of HIT_forcing

Direct numerical simulation (DNS) has been attempted since the first work reported by Orszag and Patterson \cite{OP_PRL72} for a Reynolds number based on Taylor's microscale ($\lambda$) given by $R_{\lambda} = 35$, using a $32^3$ computational periodic box. This DNS of three-dimensional (3D) homogeneous isotropic turbulence was compared with the predictions of the direct-interaction turbulence theory. With increase in computing power, $R_{\lambda}$ increased, with Kaneda and Ishihara \cite{Kaneda_Ishihara2003} reporting results for $R_{\lambda} = 1200$ using a $4096^3$ periodic box for 3D homogeneous isotropic turbulence simulations. These activities have progressed to a state where Buaria {\it et al.} \cite{Buaria} have reported results for $R_{\lambda} = 1300$ using a $12288^3$ periodic computational box for the simulation of the 3D homogeneous isotropic turbulence. All these simulations for homogeneous isotropic turbulence are performed using the incompressible Navier-Stokes equation (INSE) whose generic form is given by,

\begin{equation}
 \frac{\partial {\bf u}}{\partial t} + {\bf u}\cdot \frac{\partial {\bf u}}{\partial x} = -\nabla P/\rho + \nu \nabla^2 {\bf u} + {\bf f},
 \label{eqn1}
 \end{equation}

\noindent where ${\bf u}$ is the solenoidal velocity field, $P$ is the pressure, $\rho$ is the density of the fluid, $\nu$ is the kinematic viscosity, and ${\bf f}$ is an applied forcing imposed at large length scale. Such additional forcing for INSE has been weakly justified as necessary to ensure statistical stationarity of the computed turbulent signal by Easwaran and Pope \cite{41_Buaria} for simulation of homogeneous isotropic turbulence by the RK2-Fourier spectral method, as in Buaria {\it et al.} \cite{Buaria} This was originally used in a code by Rogallo \cite{42_Buaria} with such forcing in INSE, and the same version of the code continues to be used by a large number of researchers \cite{11_Buaria, 15_Buaria, 21_Buaria, 22_Buaria, 32_Buaria, Yueng_Donzis_KRS_JFM12, Donzis_Yueng_KRS_PoF08, Donzis_Yueng_PhysicaD10, ARanjan_PADavidson_IUTAM} among many others. This approach of using additional forcing has been critiqued \cite{TKS_PRE14} for not only for this canonical homogeneous problem, but also for so-called DNS of geophysical and engineering flows \cite{Smith_Yakhot, Bracco_McWilliams, Skote_Henningson}. For the case of transitional flows, there are many solutions of INSE which do not use such artifices and can be viewed truly as DNS \cite{Rist_Fasel_IUTAM_TKS, TKS_PRL11, TKS_PRE14}.

%@article{OP_PRL72,
%  title = {Numerical Simulation of Three-Dimensional Homogeneous Isotropic Turbulence},
%  author = {Orszag, Steven A. and Patterson, G. S.},
%  journal = {Phys. Rev. Lett.},
%  volume = {28},
%  issue = {2},
%  pages = {76--79},
%  numpages = {0},
%  year = {1972},
%  month = {Jan},
%  publisher = {American Physical Society},
%  doi = {10.1103/PhysRevLett.28.76},
%  url = {https://link.aps.org/doi/10.1103/PhysRevLett.28.76}
%}

%\bibitem{Kaneda_Ishihara2003}
%Energy dissipation rate and energy spectrum in high resolution direct numerical simulations of turbulence in a periodic box -- Yukio Kaneda and Takashi Ishihara
%Physics of Fluids 15, L21 (2003); https://doi.org/10.1063/1.1539855

%Hint of problems and finite time singularity of Euler equation

It is however noted that such simulations often display issues like blow-up at finite time (specially in the limit of vanishing $\nu$). Also "experimental measurements in homogeneous turbulence at high Reynolds numbers show the unequivocal presence of a “bottleneck” effect \cite{Lamorgese}" that has encouraged the authors in this reference to introduce hyperviscosity in the INSE as, 

\begin{equation}
 \frac{\partial {\bf u}}{\partial t} + {\bf u}\cdot \frac{\partial {\bf u}}{\partial x} = -\nabla P/\rho + (-1)^{h+1} \nu_h \nabla^{2h} {\bf u} + {\bf f},
 \label{eqn2}
 \end{equation}

\noindent where $\nu_h$ is the specified constant hyperviscosity coefficient, and ${\bf f}$ is the forcing function, so that the INSE corresponds to $h=1$ and ${\bf f}=0$. Lamorgese {\it et al.} \cite{Lamorgese} noted the "(u)se of ${\bf f}$ is unnatural (as are most other ways of forcing turbulence). However, we are primarily concerned with bottleneck effects on energy spectra, i.e., we investigate one particular characteristic of small-scale turbulence. In this case, use of a large-scale forcing (in order to analyze statistically stationary rather than decaying turbulence) is justifiable on the grounds that the details of the forcing have little effect on the small-scale statistics". 

%\bibitem{Lamorgese}
%Direct numerical simulation of homogeneous turbulence with hyperviscosity-A. G. Lamorgese, D. A. Caughey, and S. B. Pope. 
%Physics of Fluids 17, 015106 (2005); https://doi.org/10.1063/1.1833415

Bracco and McWilliams \cite{Bracco_McWilliams} solved the homogeneous stationary 2D turbulence by solving the vorticity transport equation given by,

\begin{equation}
 \frac{D \omega}{Dt}= D + {\bf F},
 \label{eqn3}
 \end{equation}

\noindent where ${\bf F}$ is again the large scale forcing term, and the other term is given by, $D = \mu \nabla^{-2} \omega + \nu \nabla^2 \omega$. The first term in this has been termed by the authors as "physically artificial hypoviscosity". These work on homogeneous, stationary turbulence in 2D and 3D requiring "hypoviscosity" and/ or "hypervisocisity" is the indication that one is not solving the true INSE. It is to be noted \cite{TKS_PRE14}, the time-averaged, compensated spectral density for (i) streamwise, (ii) wall-normal, and (iii) spanwise velocity components plotted as a function of streamwise wavenumber show the bottleneck effect noted in the experiments of Saddoughi and Veeravalli \cite{saddoughi_Veeravalli_JFM} for the inhomogeneous flow over a flat plate. In these simulations \cite{TKS_PRE14}, regular INSE has been solved in velocity-vorticity formulation by high accuracy compact schemes that requires fourth diffusion term for numerical stabilization.

There is another aspect of viscous action which is obscured when one inspects periodic flow. Authors have shown the transport equation for enstrophy ($\Omega_1 = \omega_i \omega_i$) of inhomogeneous flow to be given by \cite{TKS_RLDC}, 
$$\frac{D\Omega_1}{Dt} = \frac{2}{Re} \left[ \frac{1}{2} \nabla^2 \Omega_1 - (\nabla \omega)^2 \right]+2\omega_i \omega_j \frac{\partial u_i}{\partial x_j}$$

The first set of terms on the right hand side arises due to viscous diffusion, and furthermore, the first term ($\nabla^2 \Omega_1$) will drop out for periodic problems. This has led to a misunderstanding about turbulent flows where diffusion term is negative definite for periodic problems, as in homogeneous isotropic turbulence, leading to associating the diffusion term with dissipation, even when the flow is not periodic. The last term on the right hand side of the above equation arises due to vortex stretching, which can act as the production term of enstrophy. However, Buaria {\it et al.} \cite{Buaria} have identified this term also to additionally contribute to self-attenuation in the presence of high nonlinearity during extreme events by creating length scales smaller than the Kolmogorov scale. For Euler equation (in the limit of $Re \rightarrow \infty$ or $\nu \rightarrow 0$) the authors \cite{Buaria} have noted that enstrophy can grow unbounded in a finite time citing an observation in Doering \cite{4_Buaria}. Interestingly the authors \cite{Buaria} furthermore conjectured that such finite time blow-up "would correspond to turbulent solutions of the INSE". At the same time the authors also said that "the question remains open whether the non-linear amplification could overcome viscous damping when the flow is sufficiently turbulent". During transition to turbulence of inhomogeneous 2D flow past different bodies, authors \cite{TKS_PRE12} have shown that after flows become fully turbulent following transition, the velocity components always achieve limiting value, even when the vortex stretching is absent. 

Having dwelt upon various physical aspects of the flow, one notes that the role of numerical analysis of computing methods has received less attention. However, a very important result was brought forth for the solution of convection equation \cite{TKS_IUTAM} 

\begin{equation} 
\frac{\partial u}{\partial t} + c\frac{\partial u}{\partial x}= 0
\label{eqn4}
\end{equation}

\noindent where the Fourier spectral- RK2 method has been shown to be unconditionally unstable. This is explained next, along with the properties of the solution of the other canonical problem of convection-diffusion equation. 

The unknown in GSA is written as $u(x,t) = \int \hat{U}(k,t)e^{ikx}dk$, and if one starts solving Eq. \eqref{eqn4} with the initial condition given by $u(x_j,t = 0) = u_j^0 = \int U_0 (k)\; e^{ikx_j} dk$, then the numerical solution after $n^{th}$ time-step will be obtained as,

\begin{equation}
u_{N,j}^n = \int U_0(k)\; [|G_j|]^n\; e^{i(kx_j - n\phi_j)} dk
\label{eqn5}
\end{equation}

\noindent where $G_j = \left(\frac{\hat{U}(k,t^n+\Delta t)}{\hat{U}(k,t^n)}\right)$ is the amplification factor and is in general, a complex quantity i.e. $G_j = G_{rj} + iG_{ij}$. The modulus of this is given by $|G_j| = (G_{rj}^2 + G_{ij}^2)^{1/2}$ and the phase shift per time-step is calculated from, 
$\tan \phi_j = -{G_{ij}}/{G_{rj}}$. The GSA has been used to develop the error propagation equation for the convection equation given by \cite{TKS_vonNeumann07},

\begin{displaymath}
\frac{\partial e}{\partial t} +c \frac{\partial e}{\partial x}= -[1 -\frac{c_N}{c}] c\frac{\partial \bar{u}_N}{\partial x} - \int \frac{V_{gN} -c_N}{k} \biggl[\int ik'U_0\; [|G|]^{n}\; e^{ik' (x-c_N t)} dk'\biggr] dk \\
\end{displaymath}
\begin{equation}
\hspace{35mm} - \int \frac{{\rm Ln}~|G|}{\Delta t} U_0\; [|G|]^{n}\; e^{ik(x - c_N t)}\; dk
\label{eqn6}
\end{equation}

Here, the numerical error is given by, $e(x,t) = u(x,t) - \bar{u}_N(x,t)$, from which one can derive the governing equation for $e(x,t)$ \cite{TKS_vonNeumann07}.
The numerical phase speed ($c_N$) is obtained from $\phi_j$ (phase shift per time step), so that $n\phi_j = k c_N n \Delta t$. The physical phase speed 
is $c$ for all wavenumbers, but $c_N$ is noted to depend on $k$. Thus, the numerical solution is dispersive (in contrast to the non-dispersive nature of 1D convection equation), and is rewritten as, 

\begin{equation}
\bar{u}_N = \int U_0(k)\; [|G|]^{t/\Delta t}\; e^{ik(x - c_N t)} dk
\label{eqn7}
\end{equation}

Having obtained the correct numerical phase speed, numerical group velocity at the $j^{th}$-node is expressed as, $\biggl[ \frac{V_{gN}}{c}\biggr]_j = \frac{1}{h N_c} \frac{d\phi_j}{dk}$. Equation \eqref{eqn6} clearly demonstrates that the accuracy of the numerical solution is governed by the forcing terms on the right hand side, such as $|G|$, $c_N/c$ and $V_{gN}/c$ as the main numerical parameters for convection equation determining the error forcing. Such forcing arises due to finite discretizations, which dominates the round-off error as the other contributor. Previously, looking at the linear nature of the governing Eq. \eqref{eqn6}, researchers assumed the error dynamics also is dictated by this equation with right hand side equal to zero.  

It is important to understand the roles played by these numerical parameters in ensuring accuracy of scientific computing. The concept of error propagation was introduced in proper perspective for convection equation \cite{TKS_vonNeumann07} and the same has been shown for the diffusion equation \cite{B74} and convection-diffusion equation \cite{B94} also subsequently. GSA has also been used to demonstrate a linear dispersive mechanism for numerical error growth \cite{Tan2021}. Recently, GSA has also been used to understand the effect of dispersion and dissipation on both convection and diffusion terms in the 2D linearized compressible Navier--Stokes Equations (LCNSE) for hybrid finite difference/Fourier Spectral scheme \cite{Tan2021}.

A major development was achieved, where it was shown for the first time that the convection-diffusion equation truly reflects the dependence on numerical parameters as it does identically for INSE  \cite{B118}. This provides the framework in calibrating any numerical methods for the solution of INSE (both homogeneous and inhomogeneous flows, with or without periodicity) by checking the numerical parameters to be used in the solution process, by testing it first on the convection-diffusion equation which is given by,

\begin{equation} 
\frac{\partial u}{\partial t} + c\frac{\partial u}{\partial x}= \nu \frac{\partial^{2} u}{\partial x^{2}}
\label{eqn8}
\end{equation}  

As noted above that in numerical simulation the physical phase speed ($c$) alters to the wavenumber-dependent numerical phase speed ($c_N (k)$) for convection equation; in Eq. \eqref{eqn8}, the coefficient of diffusion (which is identical to the kinematic viscosity in INSE) also changes from constant $\nu$ to a wavenumber-dependent numerical coefficient: $\nu_N (k)$. It is customary \cite{HACM}, the numerical parameters given by the Courant-Friedrich-Lewy (CFL) number ($N_c =c\Delta t/ \Delta x$) and the Peclet number ($Pe = \nu \Delta t / (\Delta x)^2$) are important in fixing grid spacing ($\Delta x$) and time step ($\Delta t$) for solving problems dominated by convection and diffusion.

\section{Illustrative Simulations of Convection- and Convection-Diffusion Equation}
First, we discretize the governing equation using multi-stage Runge-Kutta time integration \cite{HACM} when the convection and diffusion terms are discretized by Fourier spectral method given by, $\frac{\partial u}{\partial x} = \int ik \hat{U} e^{ikx} dk$ and $\frac{\partial^2 u}{\partial x^2} = \int -k^2 \hat{U} e^{ikx} dk$. In the computations, these derivatives are evaluated using fast Fourier transform (FFT). From these discrete equations, and using the definition of numerical factors one obtains the complex parameter $G$ as function of $k\Delta x$ and $N_c$ for the convection equation \cite{TKS_vonNeumann07}; $k\Delta x$, $N_c$ and $Pe$ for convection-diffusion equation \cite{B94}. Having obtained $|G|$ and $\phi$, one can readily obtain $c_N/c$ and $V_{g,N}/c$ for the convection equation. For the convection-diffusion equation, one would also find out $\nu_N/\nu$, the non-dimensional dispersive coefficient of diffusion \cite{B94, B118}.

First, we demonstrate solution of convection equation obtained by solving Eq. \eqref{eqn6} by Fourier spectral method for spatial discretization and use two-, three-, and four-stage Runge-Kutta method. If one chooses a grid spacing of $\Delta x$, then the resolved maximum wavenumber ($k_{max}$) is fixed by the Nyquist criterion, 
$k_{max} \Delta x = \pi$. The convection equation is solved in a domain, $0 \le x \le 10$ discretized uniformly with 4096 points for the propagation of a wave-packet. The initial condition for the wave packet is given by,

\begin{equation}
u(x,0) = e^{-10(x-5)^2} \sin (k_0 \Delta x)
\label{eqn9}
\end{equation}

The wave-packet is characterized by the central wavenumber given by, $k_0 \Delta x = 0.22$, a value which is significantly lower than the Nyquist limit. This low value ensures that the error sources due to phase error and dispersion error are significantly sub-dominant as compared to the error that is created due to the term 
$|G|$. 

\begin{figure}[H]
        \centering
         \includegraphics[scale=0.5]{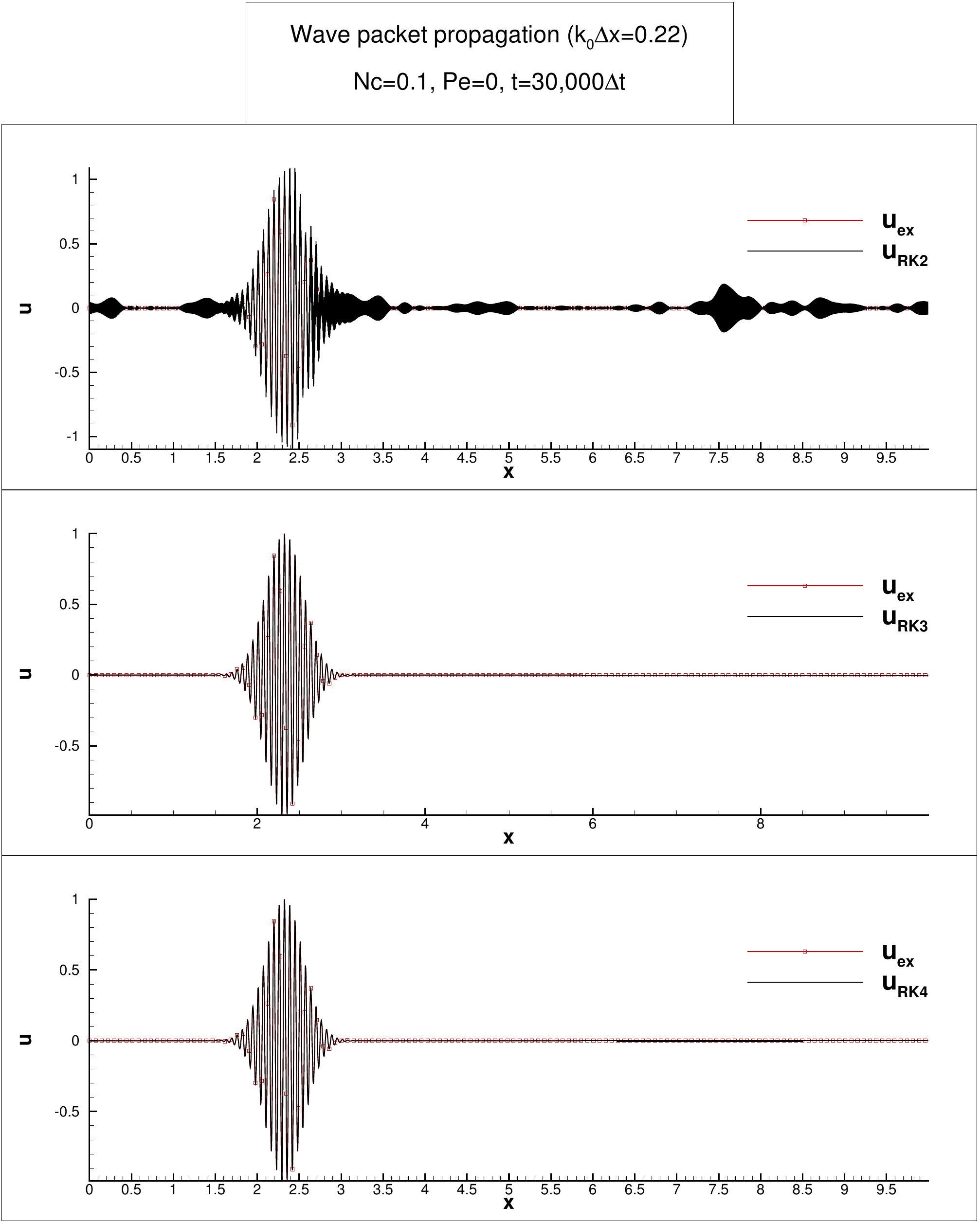}
        \caption{Comparison of numerical solution with exact solution for a time $t=30,000\Delta t$ for the convection of a wavepacket using a 1D linear convection equation. Here, Fourier spectral method is adopted for spatial derivative and RK2, RK3 and RK4 methods are employed for time integration. All simulations are performed on a uniform grid with $4096$ points and for a CFL number $N_c=0.1$.}
        \label{fig1}
    \end{figure}

In Fig. \ref{fig1}, the propagation of the wave-packet is shown after 30000 $\Delta t$, where $\Delta t = 4.884 \times 10^{-4}$ is fixed from $N_c=0.1$ and $c = 0.5$, with RK2 method result shown in the top frame; the three-stage RK (RK3) method result in the middle frame and four-stage RK (RK4)  method result depicted in the bottom frame. As the exact solution is simply translated to the right at the speed $c$, one can compare it with the numerical solution. It is evident that RK2-Fourier spectral method results are erroneous. As one is solving a periodic problem, the round-off error magnifies, as determined by $G$, and that keeps convecting within the domain accumulating in magnitude. However, no distinguishable errors are noted for the RK3-Fourier spectral and RK4-Fourier spectral methods in Fig. \ref{fig1} for this low value of CFL number, $N_c= 0.1$. 

\begin{figure}[H]
        \centering
         \includegraphics[scale=0.5]{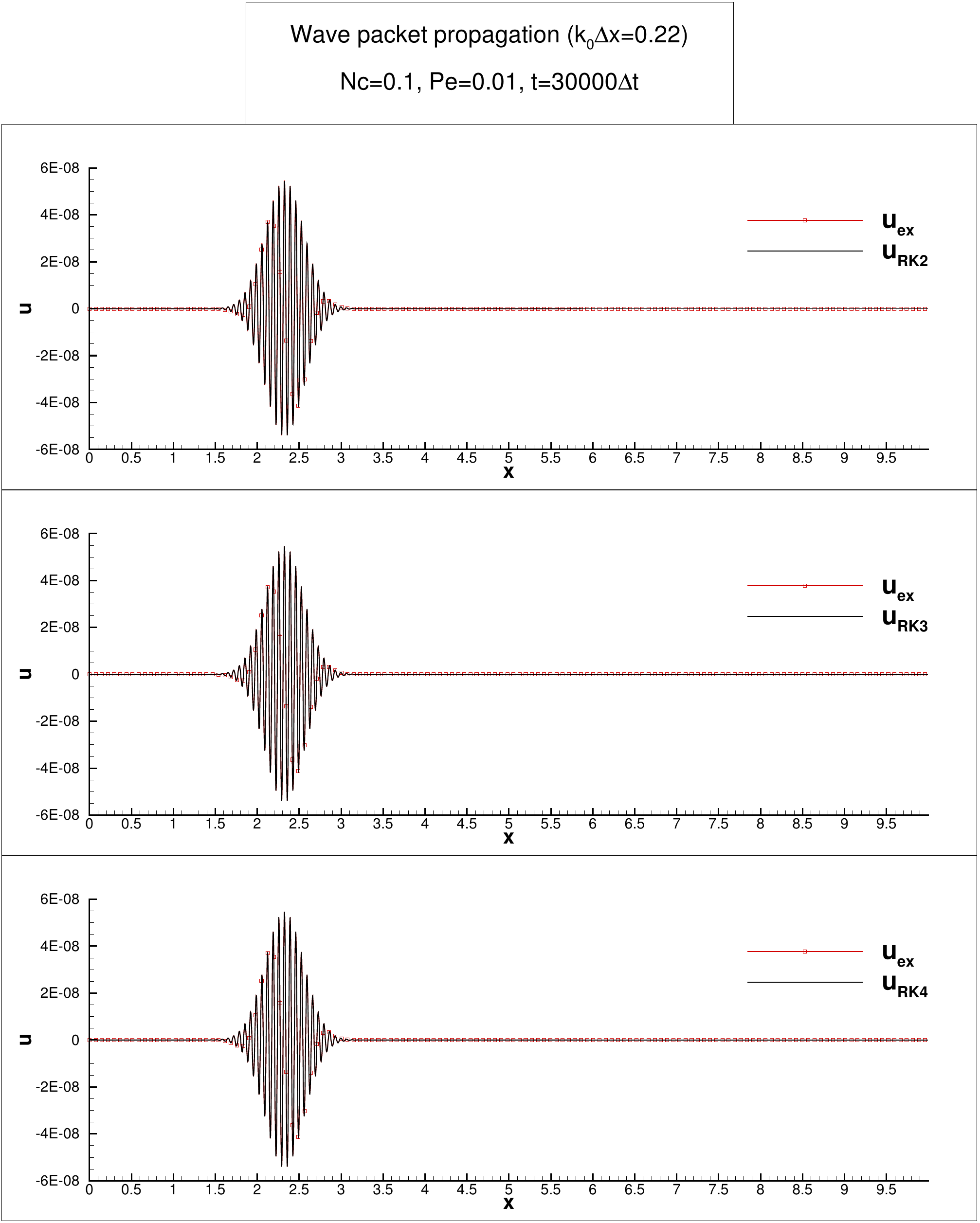}
        \caption{Comparison of numerical solution with exact solution for a time $t=30,000\Delta t$ for the convection of a wavepacket using a 1D linear convection-diffusion equation. Here, Fourier spectral method is adopted for spatial derivative and RK2, RK3 and RK4 methods are employed for time integration. All simulations are performed on a uniform grid with $4096$ points, with a CFL number $N_c=0.1$ and Peclet number $Pe=0.01$.}
        \label{fig2}
    \end{figure}

In Fig. \ref{fig2}, the propagation of the same wave-packet is shown after a time interval of 30000 $\Delta t$, for the same value of $N_c = 0.1$ with the three time discretization methods used in Fig. \eqref{fig1}, with Fourier spectral method for spatial discretization, by solving the convection-diffusion equation. Here, the same time step is used, and the coefficient of diffusion is so chosen that the Peclet number is 0.01. Unlike the non-dissipative, non-dispersive convection equation case in Fig. \ref{fig1} for which the physical amplification rate is given by $|G_{phys}| = 1$, here for the convection-diffusion equation case the physical amplification would indicate attenuation due to diffusion, such that  $|G_{Phys}| < 1$. In this case, the numerical amplification rate ($|G_{Num}|$) will be such that $|G_{Num}|/|G_{Phys}|$ is equal to one for accuracy, and $|G_{Num}| > 1$ would indicate numerical instability. In Fig. \ref{fig2}, the solutions are shown for $N_c = 0.1$ and $Pe = 0.01$ for the three Runge-Kutta methods, with the same wave-packet starting from the same initial locations. All the numerical solutions coincide with the exact solution, implying that the unstable RK2-Fourier spectral method for convection equation has become stable for the convection-diffusion equation for the same value of $N_c = 0.1$. It has already been shown unambiguously that the numerical behavior of convection-diffusion equation reflects identically the behavior of the numerical solution of INSE when the numerical parameters are kept identical \cite{B118}. For the case of RK3 and RK4 methods, the solutions are error free for these parameter combinations for both the convection- and convection-diffusion equations.

    %\clearpage
    
    %\clearpage
    \begin{figure}[H]
        \centering
         \includegraphics[scale=0.5]{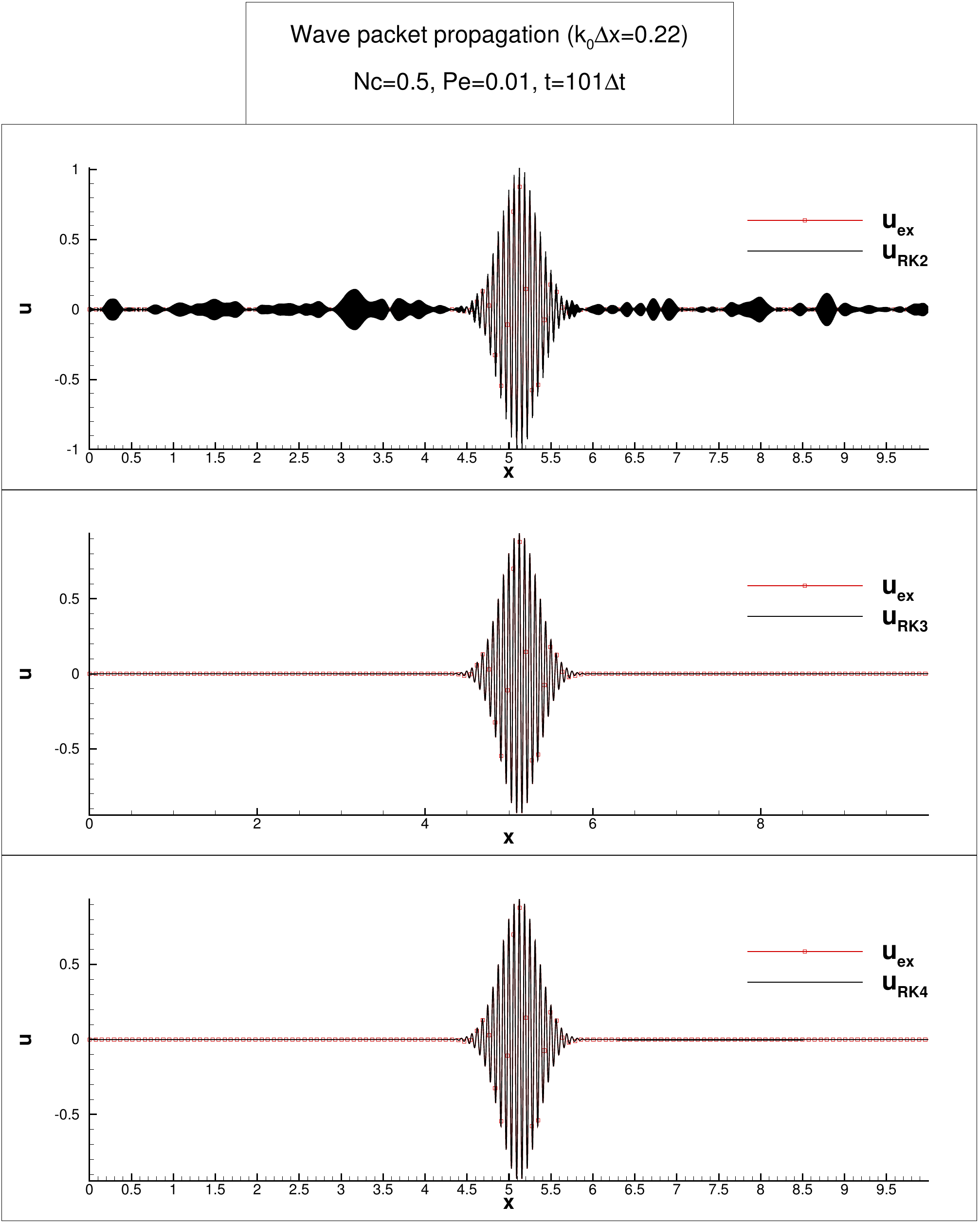}
        \caption{Comparison of numerical solution with exact solution for a time $t=101\Delta t$ for the convection of a wavepacket using a 1D linear convection-diffusion equation. Here, Fourier spectral method is adopted for spatial derivative and RK2, RK3 and RK4 methods are employed for time integration. All simulations are performed on a uniform grid with $4096$ points, with a CFL number $N_c=0.5$ and Peclet number $Pe=0.01$.}
        \label{fig3}
    \end{figure}
    
The conditional stability of the RK2-Fourier spectral method for the solution of convection-diffusion equation, explains its use for the simulation of INSE, despite its failure to do so for the convection equation. These aspects are further explored by solving convection-diffusion equation for a higher CFL number of $N_c = 0.5$, while keeping the Peclet number at the same level of $Pe = 0.01$. It is well known \cite{20_Buaria} that for DNS, the diffusion discretization places a stricter control on the time step via the choice of Peclet number. The solutions obtained by RK2, RK3 and RK4 methods are compared in Fig. \ref{fig3} for $N_c= 0.5$ and $Pe = 0.01$ and the simulations are run only for 101 time steps. One notices erroneous solutions by the RK2-Fourier spectral method, whereas RK3- and RK4-Fourier spectral methods again reproduce error free propagation of the same wave-packet as before. All of these results shown here indicate the necessity to characterize the space-time discretization methods. Only characterization in this respect has been reported in Fig. 1 by Sengupta \cite{TKS_IUTAM} for the RK2-Fourier spectral method for the solution of convection equation. It is shown that this method (which is used by numerous researchers and continues to be used for DNS of INSE), is unconditionally unstable for all CFL numbers for flows highly dominated by convection. Only for very small values of $k\Delta x$ one observes neutral stability, which can explain as to why researchers \cite{5_Buaria} have reported finite-time blow-up for the solution of Euler equation. To understand the implications of the reported results in Figs. \ref{fig1} to \ref{fig3} one must perform a thorough study of the numerical properties of convection and convection-diffusion equations for the pseudo-spectral method proposed by Rogallo \cite{42_Buaria} using RK2 time integration. The numerical evidences provided in Figs. \ref{fig1} to \ref{fig3} indicate that the shortcoming of the RK2-Fourier spectral method can be possibly overcome by using RK3 and/ or RK4 time integration methods.

\section{Error Dynamics for Convection and Convection-Diffusion Equations}

It has been already noted in Eq. \eqref{eqn6} that there are three main sources driving the dynamics of error, via the forcing by the phase error term, the dispersion error term and more importantly by numerical instability and unphysical numerical stability term, noted in Eq. \eqref{eqn6} for the convection equation. It is a common misconception that numerical stability is desirable \cite{vonNeumann_Richtmeyer}, which makes many practitioners complacent about this potential source of error, whereas a numerically unstable method is easily recognized due to solution blow-up. For the present analysis of convection and convection-diffusion equation, numerical parameters are chosen such that the error is solely due to numerical instability and unphysical numerical stability. For the convection equation, an error-free calculation must have neutrally stable numerical amplification factor, $|G|=1$. In contrast, for convection-diffusion equation the physical solution attenuates such that $|G_{Phys}| < 1$ and the numerical solution should reflect this, so that $|G_{Num}|/|G_{Phys}| = 1$. The procedures to calculate these quantities are given for convection equation \cite{TKS_vonNeumann07} and convection-diffusion equation \cite{B94}.  

    \begin{figure}[H]
        \centering
         \includegraphics[scale=0.7]{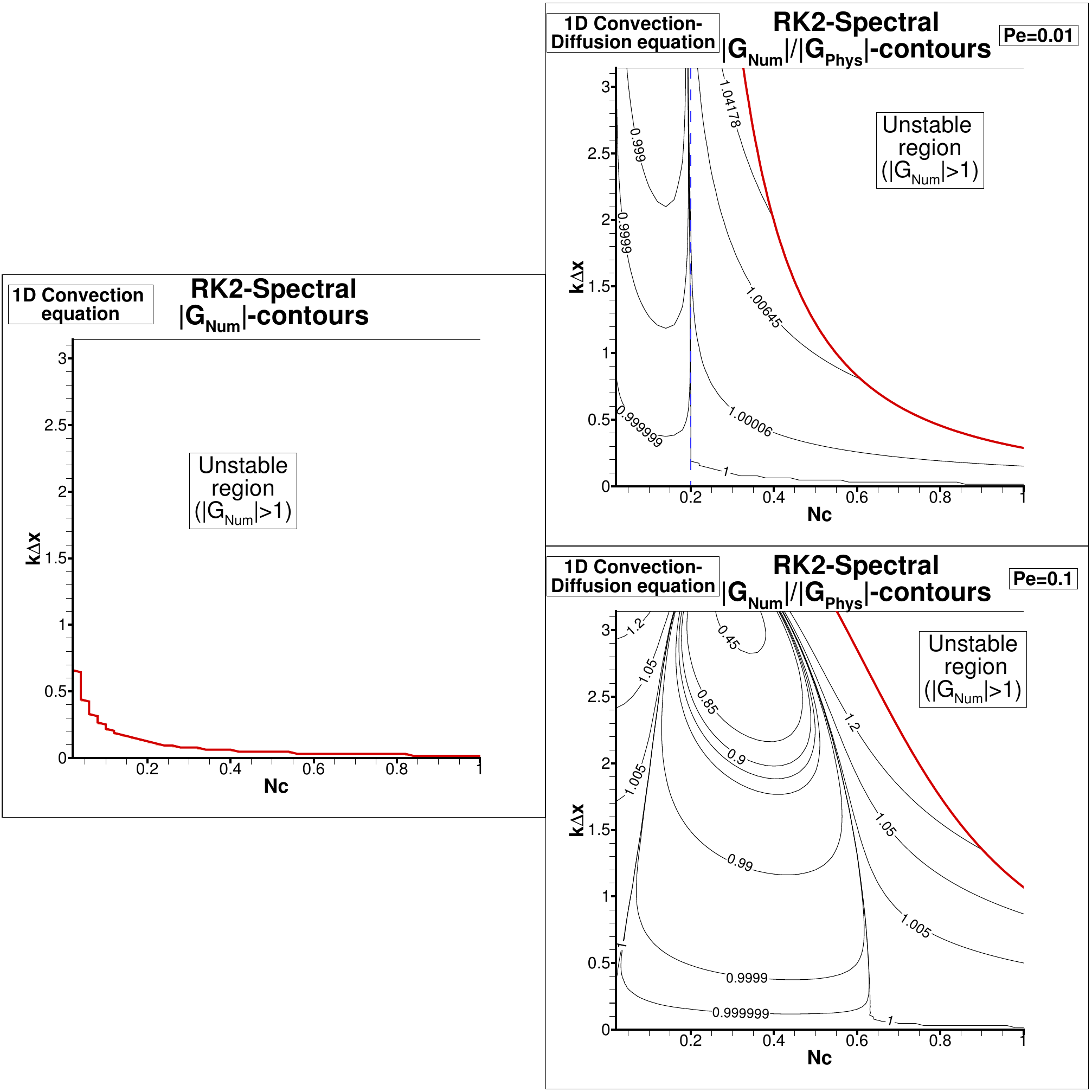}
        \caption{Numerical amplification factor and ratio of numerical amplification factor to the physical amplification factor for the RK2-Fourier spectral scheme plotted in the $(N_c,k\Delta x)$-plane for linear convection equation and linear convection-diffusion equation, respectively. For the latter equation, contours are plotted for two representative Peclet numbers- $Pe=0.01$ and $0.1$. Dashed line represents the optimal $N_c$ value.}
        \label{fig4}
    \end{figure}
    
We use the GSA to map the property charts of the method involving Fourier spectral method for spatial discretization and RK2 method for time integration. We also report the properties at the interior of the domain, as one is interested to analyze the methods to solve 3D homogeneous isotropic turbulence in a periodic box. In Fig. \ref{fig4}, RK2-Fourier spectral method is considered for both the canonical equations. For the one-dimensional convection equation, $|G|$-contours are shown plotted in the $(N_c, k\Delta x)$-plane on the left side of Fig. \ref{fig4}. It is noted that barring a small region very close to the origin, the method is unconditionally unstable. This result was also reported \cite{TKS_IUTAM}, and one notes the effects of numerical instability in Fig. \ref{fig1}, in the top frame for which the choice of $N_c =0.1$ and $k_0 \Delta x = 0.22$ does not guarantee stability of all the resolved scales by the Nyquist criterion. The seeds of unstable signals are always present in the round-off error and those wavenumber components interact with the given wave-packet to create multitudes of wave-packets noted for this periodic problem. A non-periodic problem with different inflow and outflow boundary conditions would allow some of the error-packets to leave the computational domain. To emphasize the growth of error that is inherent with the numerical method, a periodic problem is considered so that the signal is always inside the domain, along with the recirculating continuously evanescent error. For the convection-diffusion equation, the error dynamics is additionally affected via a term involving $\nu_N/\nu$ and the error dynamics is given by,

 \begin{equation}
\begin{split}
e_t + ce_x - \nu e_{xx} =& \int_{-k_{max}}^{k_{max}} (\nu_N-\nu) \; k^2 \; e^{-\nu_N k^2 n \Delta t} \; U_0(k) \; e^{ik(x-c_Nt^n)}\,dk \\
& +ikc_N \; e^{-\nu_N k^2 n \Delta t} \; U_0(k) \; e^{ik(x-c_Nt^n)} \biggl. \biggr|^{k_{max}}_{-k_{max}} \\
&- \int_{-k_{max}}^{k_{max}} \left( \frac{V_{gN}-c_N}{k} \right) \left\{ \int_{-k_{max}}^{k} ik^{\prime} \; e^{-\nu_N {k^{\prime}}^2 n \Delta t} \; U_0(k^{\prime}) \; e^{ik^{\prime}(x-c_Nt^n)} \;dk^{\prime} \right\} dk\\
&- \int_{-k_{max}}^{k_{max}} ikc \; e^{-\nu_N k^2 n \Delta t} \; U_0(k) \; e^{ik(x-c_Nt^n)} \;dk
\end{split}
\label{eqn10}
\end{equation}

On the right hand side frames of Fig. \ref{fig4}, we have shown the contours of $|G_{Num}|/|G_{Phys}|$ in ($N_c, k\Delta x$)-plane for $Pe = 0.01$ and 0.1. One notices that there is a finite range of $N_c$ for which there is no instability over the full resolved scale. The instability is indicated whenever $|G_{Num}| > 1$. This is true for both the Peclet numbers, but the contour values in the stable region varies with varying $k\Delta x$. Presented results here are for $N_c= 0.1$ and 
0.5, for which the quotient is almost equal to $1$ for only a small range of $k\Delta x$ within the Nyquist limit. To find such a value of $N_c$ for a chosen $Pe$ with maximum range of wavenumbers within the Nyquist limit, one requires to perform GSA. The authors are not aware that such an activity has been performed before by any researchers.

    \begin{figure}[H]
        \centering
         \includegraphics[scale=0.7]{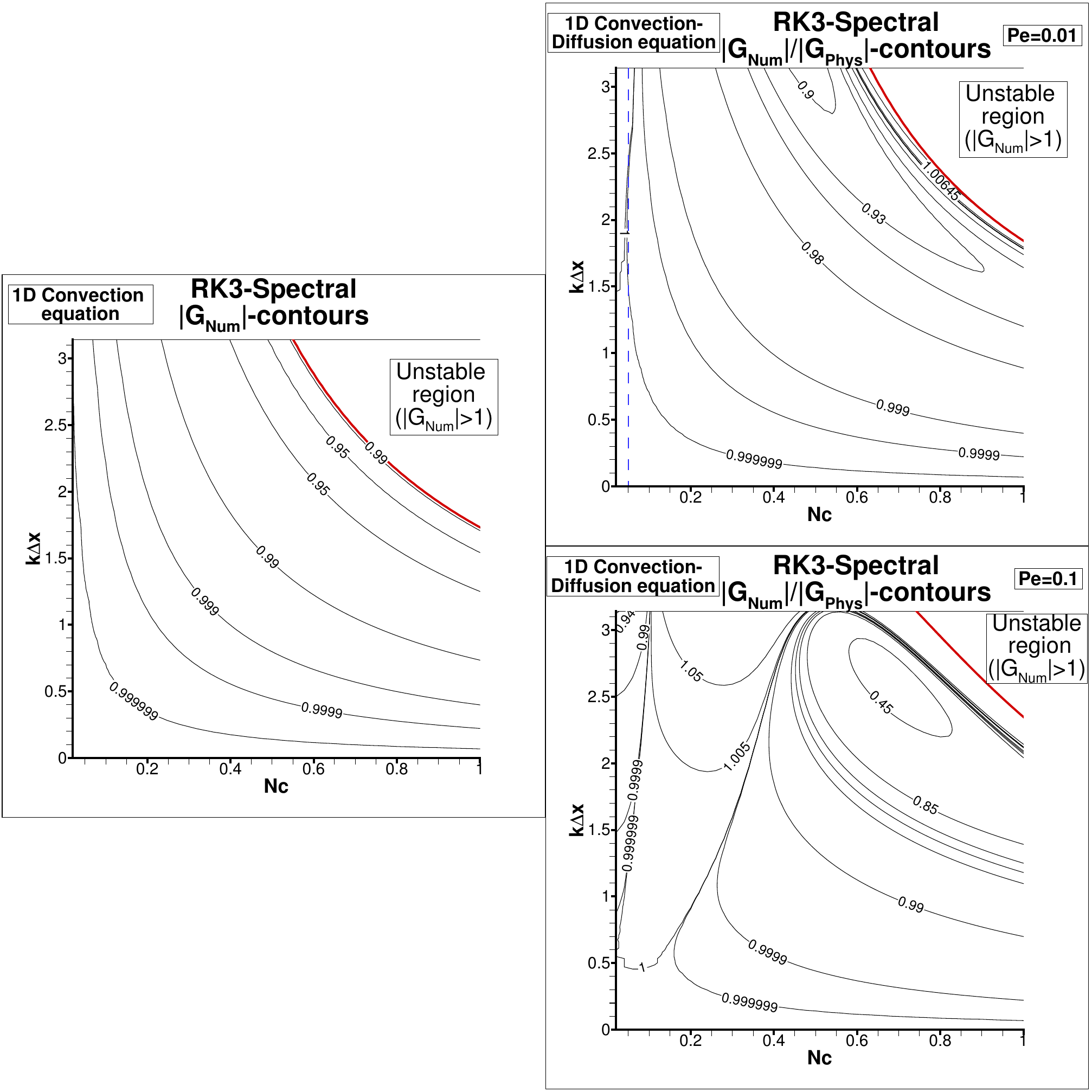}
        \caption{Numerical amplification factor and ratio of numerical amplification factor to the physical amplification factor for the RK3-Fourier spectral scheme plotted in the $(N_c,k\Delta x)$-plane for linear convection equation and linear convection-diffusion equation, respectively. For the latter equation, contours are plotted for two representative Peclet numbers- $Pe=0.01$ and $0.1$. Dashed line represents the optimal $N_c$ value for LES.}
        \label{fig5}
    \end{figure}
    
        \begin{figure}[H]
        \centering
         \includegraphics[scale=0.7]{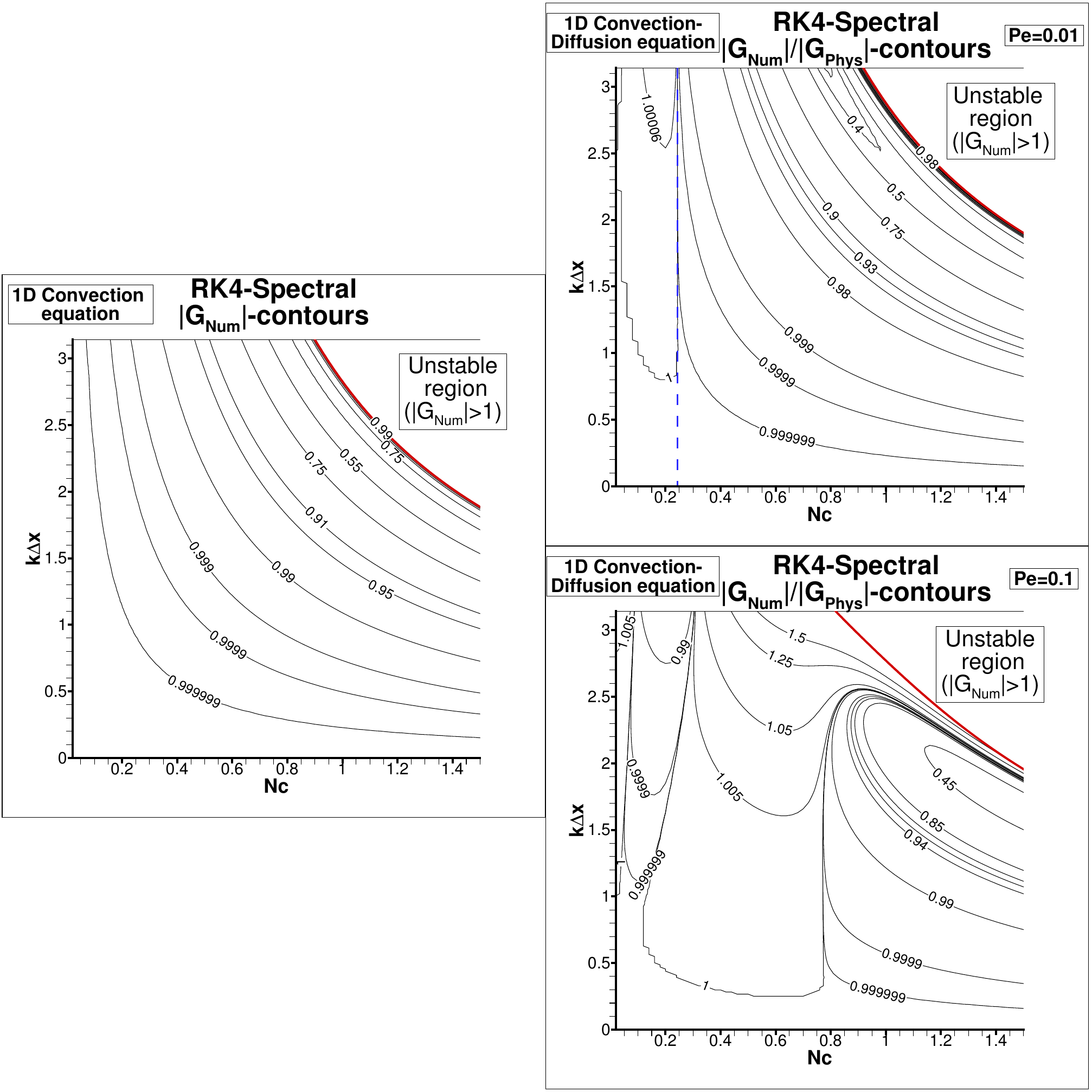}
        \caption{Numerical amplification factor and ratio of numerical amplification factor to the physical amplification factor for the RK4-Fourier spectral scheme plotted in the $(N_c,k\Delta x)$-plane for linear convection equation and linear convection-diffusion equation, respectively. For the latter equation, contours are plotted for two representative Peclet numbers- $Pe=0.01$ and $0.1$. Dashed line represents the optimal $N_c$ value for LES.}
        \label{fig6}
    \end{figure}
    
In Figs. \ref{fig5} and \ref{fig6}, $|G_{Num}|/|G_{Phys}|$-contours are plotted in ($N_c, k\Delta x$)-plane for $Pe = 0.0$, 0.01 and 0.1 for RK3- and RK4-Fourier spectral methods. Unlike the RK2-Fourier spectral method for the convection equation, RK3 and RK4 time integration methods do not exhibit unconditional numerical instability, as shown in the left panel. However to get neutral stability for this canonical equation, one will be forced to take small value of $N_c$, with RK4 method significantly better than the RK3 time integration method allowing much higher values of $N_c$. Such values of CFL number will not display finite time blow-up of solution for Euler equation. The right frames in Figs. \ref{fig5} and \ref{fig6}, $|G_{Num}|/|G_{Phys}|$-contours are plotted in ($N_c, k\Delta x$)-plane for $Pe = 0.01$ and 0.1 for simulating INSE. For RK3 time integration method, performing accurate simulation would require taking vanishingly small values of $N_c$, for vanishingly small values of $Pe$. It must however, be remembered that for no parameter combinations, one will reproduce the exact solution for the fully resolved Nyquist limit of wavenumbers 
($k\Delta x = \pi$). For the RK4 time integration method, the permissible values of $N_c$, $Pe$ for performing near-DNS simulation will be relatively higher than that will be obtained by the RK3 time integration method. Such imperfectly resolved computations can be viewed more as an implicit large eddy simulation. Overall, comparing the $|G_{Num}|/|G_{Phys}|$-contours in Figs. \ref{fig4} to \ref{fig6}, it appears that RK4-Fourier spectral method will provide the best results for $Pe = 0.01$ and $N_c= 0.243793$ as described next. 

\begin{figure}[H]
        \centering
         \includegraphics[width=\textwidth,keepaspectratio=true]{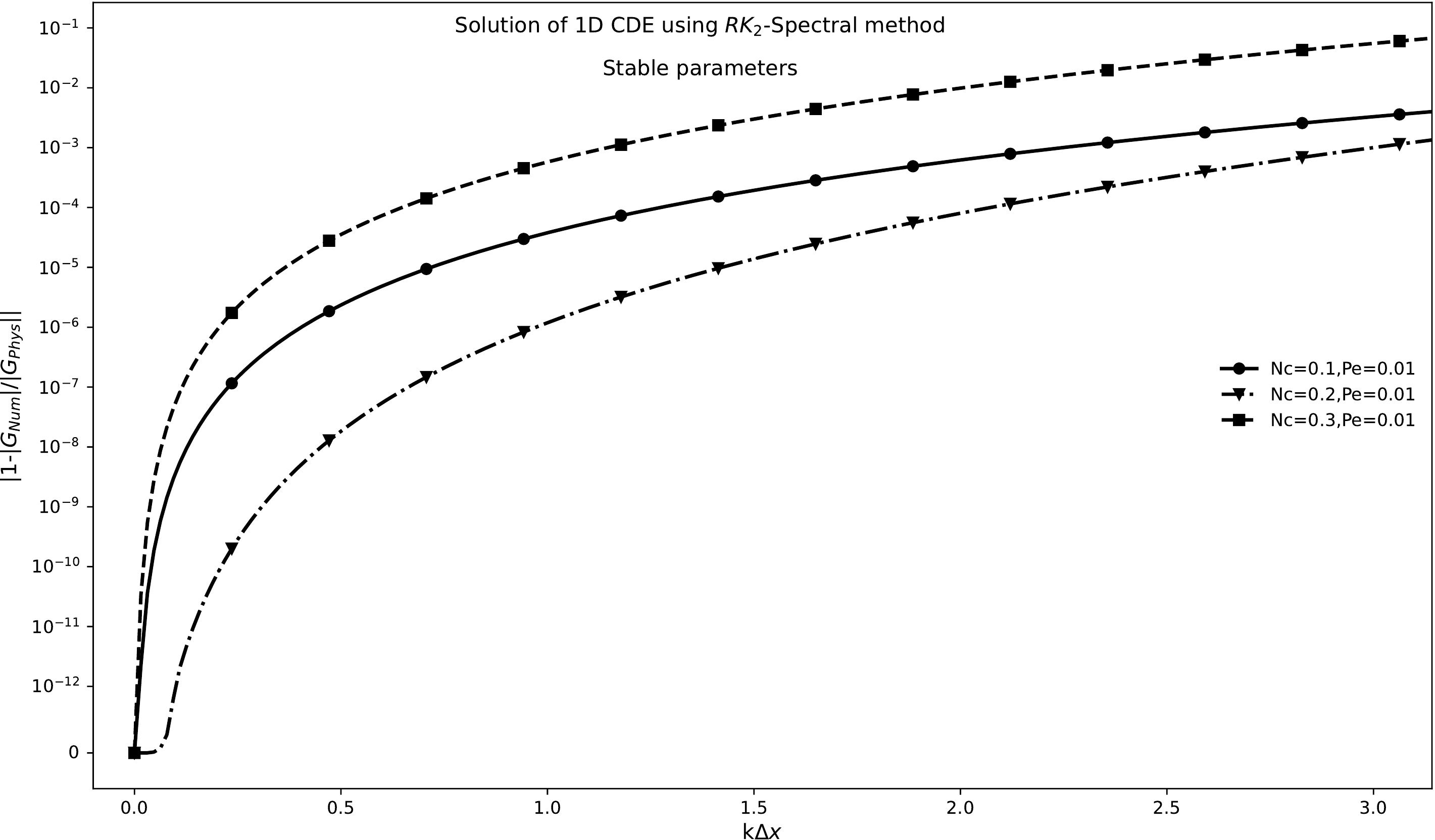}
        \caption{Determination of optimal CFL number $N_c$ for accurate solution of 1D linear convection-diffusion equation using RK2-Fourier spectral method for Peclet number $0.01$. Accuracy is evaluated with respect to error in numerical amplification factor ($|1-(|G_{num}|/|G_{phys}|)|$) for all resolvable wavenumbers $(k\Delta x)$.}
        \label{fig7}
    \end{figure}

In Figs. \ref{fig4} to \ref{fig6}, one can note the contribution to error dynamics by the term involving $|G_{Num}|/|G_{Phys}|$ for specific values of $Pe$. On a cursory observation, one notes that for a specific value of $N_c$, $|G_{Num}|/|G_{Phys}|$ value keeps changing with $k\Delta x$. To quantify such variations and to obtain an optimal value of $N_c$ for which the error will be least due to this term $|G_{Num}|/|G_{Phys}|$, in Fig. \ref{fig7}, we have plotted $||G_{Num}|/|G_{Phys}|-1|$ as a function of $k\Delta x$ for the RK2-spectral method showing that the optimal error forcing is obtained for $N_c=0.2$, which is shown to be lower than that for $N_c= 0.1$ and $0.3$. Having noted this optimal value of $N_c=0.2$, one also notices that the forcing by this term increases monotonically in nonlinear manner with increase in wavenumber up to the Nyquist limit. Such monotonic growth in error with wavenumber is clearly apparent from Fig. \ref{fig4}, for these lower values of CFL, which is also restricted to the right due to numerical instability for the choice of $N_c$. Even though the numerical instability relaxes somewhat with increase in $Pe$, the increased error with increase in wavenumber appears to be much more drastic, and should be completely avoided.

\begin{figure}[H]
        \centering
         \includegraphics[width=\textwidth,keepaspectratio=true]{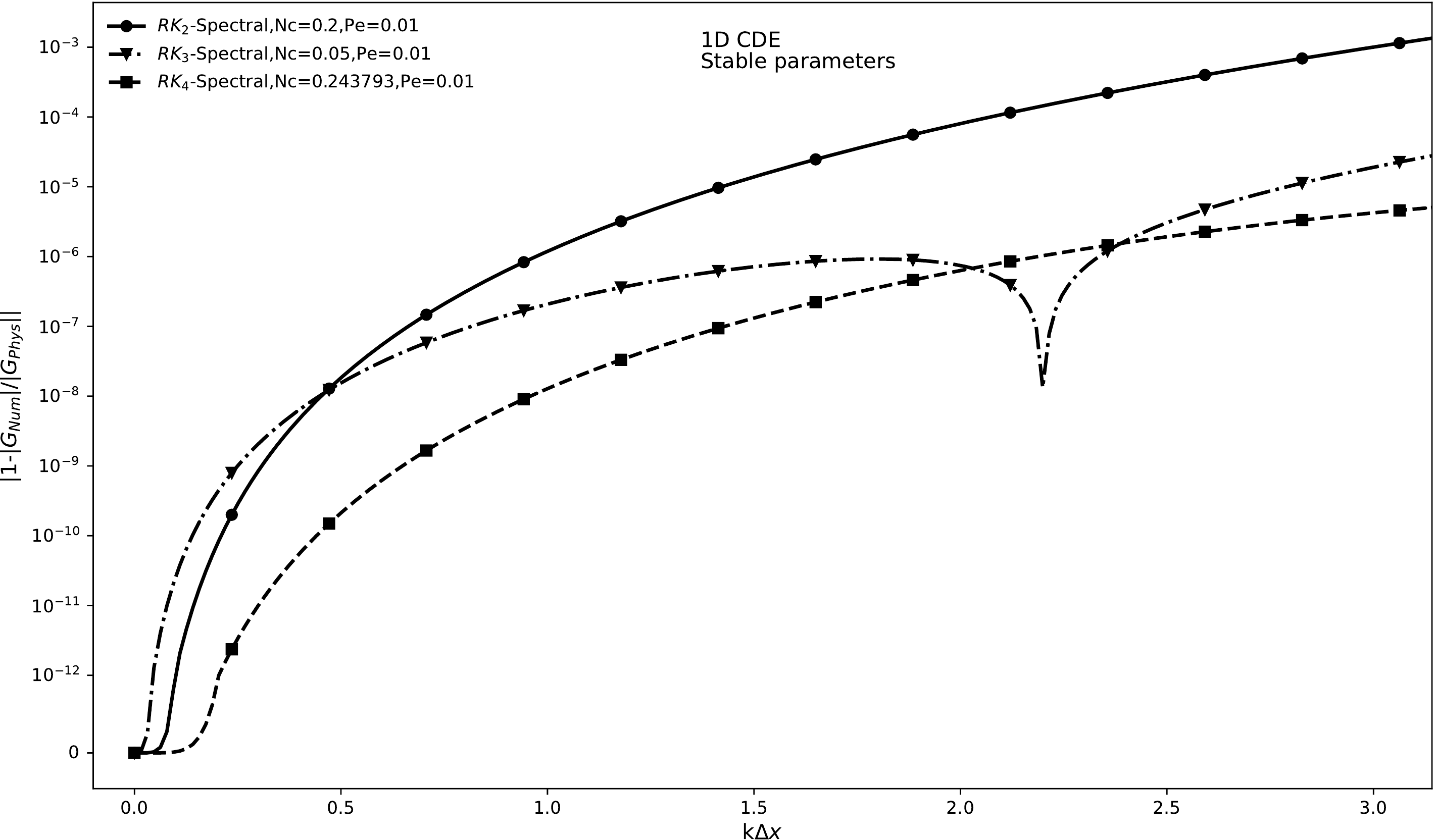}
        \caption{Comparison of the accuracy of RK2, RK3 and RK4 time integration schemes in solving the 1D linear convection-diffusion equation for Peclet number $0.01$ using Fourier spectral method. Accuracy is evaluated with respect to error in numerical amplification factor ($|1-(|G_{num}|/|G_{phys}|)|$) for all resolvable wavenumbers $(k\Delta x)$. Presented values of $N_c$ correspond to the respective optimal values of the time integration schemes.}
        \label{fig8}
    \end{figure}

Having noted the existence of optimal $N_c$ value for RK2-Fourier spectral method as equal to 0.2, next we locate the same for RK3-Fourier spectral and RK4-Fourier spectral methods. The variation of contributions by the optimal $N_c$ choices to error dynamics is compared in Fig. \ref{fig8}. In Fig. \ref{fig5}, the variation of $|G_{Num}|/|G_{Phys}|$-contours in $(N_c,k\Delta x)$-plane indicates that the error dynamics term dependent on $|G_{Num}|/|G_{Phys}|$ is more complex and contribute higher values. For the RK3-Fourier spectral method such an optimal value is noted for $N_c=0.05$ and for RK4-Fourier spectral method the optimal error is obtained for $N_c= 0.243793$ for the choice of $Pe =0.01$ with the results depicted in Fig. \ref{fig8}. The variation of error forcing for RK3-spectral method in this figure quantitatively is also apparent qualitatively in Fig. \ref{fig5}. But, the RK4-spectral has monotonic variation of $|G_{Num}|/|G_{Phys}|$ with wavenumber and that causes the least error among these three time integration methods used with Fourier spectral method for spatial discretization. 
    
\section{Summary and Conclusions}

%Suman add the points which one can take home from this study! That should be the summary and conclusion. 
A detailed analysis and calibration is performed for the first time using global spectral analysis (GSA) for the numerical methods employing Fourier spectral discretization for spatial derivatives and RK2, RK3 and RK4 methods for temporal discretization in solving the linear convection and linear convection-diffusion equations. Numerical analysis confirms the unconditionally unstable nature of RK2-Fourier spectral method in solving the 1D convection equation, as reported earlier \cite{TKS_IUTAM}. It is also demonstrated that RK3/RK4-Fourier spectral methods are stable in solving the same problem. The solution of convection-diffusion equation is interesting, as this model equation reflects truly upon the behaviour of incompressible Navier-Stokes equation. The RK2-Fourier spectral method for convection-diffusion equation does not exhibit unconditional numerical instability, as noted for convection equation. For a Peclet number equal to $0.01$ and for low values of CFL number ($N_c$), one notices numerically stable property. Error dynamics analysis is performed to answer the question of which among RK2, RK3 and RK4 methods offer the best accuracy for the pseudo-spectral method in solving the convection-diffusion equation. By noting the error contributed by numerical amplification factor with respect to physical amplification factor $|1-|G_{Num}|/|G_{Phys}||$ at all resolvable scales, optimal $N_c$ values are determined for a fixed Peclet number for the RK2, RK3 and RK4 time integration methods. Comparing such an error component due to amplification factor at the optimal $N_c$, the RK4 method shows superior accuracy. The optimal values are obtained as the following. For RK2-Fourier spectral method $N_c = 0.2$ for $Pe =0.01$. For the RK3-Fourier spectral method the optimal numerical parameter is $N_c = 0.05$ for $Pe= 0.01$ and for RK4-Fourier spectral method, these are given by, $N_c= 0.243793$ for $Pe= 0.01$. The analysis suggests that one can use Fourier spectral method by adopting RK4 method for time integration in order to obtain good accuracy. In all such cases, the full range of resolved wavenumbers given by the Nyquist limit of $k_{max} = \pi/ \Delta x$ will not be available, and calling even the optimal parameter cases will be more apt as implicit large eddy simulation, specially while tracking extreme events in turbulent flow simulations.

\section*{Data Availability}
The data that support the findings of this study are available from the corresponding author upon reasonable request.


\nocite{*}

\begin{thebibliography}{99.}

\bibitem{OP_PRL72}
Orszag S., \& Patterson G., Numerical Simulation of Three-Dimensional Homogeneous Isotropic Turbulencer, {\em Phys. Rev. Lett.}, \textbf{28}(2), 76-79 (1972)

\bibitem{Kaneda_Ishihara2003}
Ishihara T., Morishita K., Yokokawa M., Uno A. \& Kaneda Y., Energy spectrum in high-resolution direct numerical simulations of turbulence, {\em Phys. Rev. Fluids}, \textbf{1}(8), 082403 (2016)

\bibitem{Buaria}
Buaria D., Pumir A. \& Bodenschatz E., Self-attenuation of extreme events in Navier--Stokes turbulence, {\em Nat. Commun.}, \textbf{11}, 5852 (2020)

\bibitem{41_Buaria}
Eswaran V., \& Pope S. B., An examination of forcing in direct numerical simulations of turbulence, {\em Computers \& Fluids}, \textbf{16}(3), 257-278 (1988)

\bibitem{42_Buaria}
Rogallo R. S., Numerical experiments in homogeneous turbulence, NASA Tech. Memo, 81315 (1981)

\bibitem{11_Buaria}
Ishihara T., Gotoh T., \& Kaneda Y., Study of High--Reynolds Number Isotropic Turbulence by Direct Numerical Simulation, {\em Ann. Rev. Fluid Mech.}, \textbf{41}, 165-180 (2009)

\bibitem{15_Buaria}
Buaria D., Bodenschatz E., \& Pumir A., Vortex stretching and enstrophy production in high Reynolds number turbulence, {\em Phys. Rev. Fluids}, \textbf{5}(10), 104602 (2020)

\bibitem{21_Buaria}
Kaneda Y., Ishihara T., Yokokawa M., Itakura K., \& Uno A., Energy dissipation rate and energy spectrum in high resolution direct numerical simulations of turbulence in a periodic box, {\em Phys. Fluids}, \textbf{15}(2), L21-L24 (2003)

\bibitem{22_Buaria}
Buaria D., \& Sreenivasan K. R., Dissipation range of the energy spectrum in high Reynolds number turbulence, {\em Phys. Rev. Fluids}, \textbf{5}(9), 092601 (2020)

\bibitem{32_Buaria}
Choi Y., Kim B., \& Lee C., Alignment of velocity and vorticity and the intermittent distribution of helicity in isotropic turbulence, {\em Phys. Rev. E}. \textbf{80}(1), 017301 (2009)

\bibitem{Yueng_Donzis_KRS_JFM12}
Yeung P. K., Donzis D. A., \& Sreenivasan K. R., Dissipation, enstrophy and pressure statistics in turbulence simulations at high Reynolds numbers, {\em J. Fluid Mech.}, \textbf{700}, 5-15 (2012)

\bibitem{Donzis_Yueng_KRS_PoF08}
Donzis D. A., Yeung P. K., \& Sreenivasan K. R., Dissipation and enstrophy in isotropic turbulence: Resolution effects and scaling in direct numerical simulations, {\em Phys. Fluids}. \textbf{20}(4), 045108 (2008)

\bibitem{Donzis_Yueng_PhysicaD10}
Donzis D. A., \& Yeung P. K., Resolution effects and scaling in numerical simulations of passive scalar mixing in turbulence, {\em Physica D: Nonlinear Phenomena}, \textbf{239}(14), 1278-1287 (2010)

\bibitem{ARanjan_PADavidson_IUTAM}
Ranjan A., \& Davidson P. A., DNS of a buoyant turbulent cloud under rapid rotation, {in:\em IUTAM Symp. Proc. Advances In Computation, Modeling And Control Of Transitional And Turbulent Flows}, 452-460 (2016)

\bibitem{TKS_PRE14}
Bhaumik S., \& Sengupta T. K., Precursor of transition to turbulence: Spatiotemporal wave front, {\em Phys. Rev. E}, \textbf{89}(4), 043018 (2014)

\bibitem{Smith_Yakhot}
Smith L. M., \& Yakhot V., Finite-size effects in forced two-dimensional turbulence, {\em J. Fluid Mech.}, \textbf{274}, 115-138 (1994)

\bibitem{Bracco_McWilliams}
Bracco A., \& McWilliams J. C., Reynolds-number dependency in homogeneous, stationary two-dimensional turbulence, {\em J. Fluid Mech.}, \textbf{646}, 517-526 (2010)

\bibitem{Skote_Henningson}
Skote M., \& Henningson D. S., Direct numerical simulation of a separated turbulent boundary layer, {\em J. Fluid Mech.}, \textbf{471}, 107-136 (2002)

\bibitem{Rist_Fasel_IUTAM_TKS}
Rist U., \& Fasel H., Direct numerical simulation of controlled transition in a flat-plate boundary layer, {\em J. Fluid Mech.}, \textbf{298}, 211-248 (1995)

\bibitem{TKS_PRL11}
Sengupta T. K., \& Bhaumik S., Onset of Turbulence from the Receptivity Stage of Fluid Flows, {\em Phys. Rev. Lett.}, \textbf{107}(15), 154501 (2011)

\bibitem{Lamorgese}
Lamorgese A. G., Caughey D. A., \& Pope S. B., Direct numerical simulation of homogeneous turbulence with hyperviscosity, {\em Phys. Fluids}, \textbf{17}(1), 015106 (2005)

\bibitem{saddoughi_Veeravalli_JFM}
Saddoughi S. G., \& Veeravalli S. V., Local isotropy in turbulent boundary layers at high Reynolds number, {\em J. Fluid Mech.}, \textbf{268}, 333-372 (1994)

\bibitem{TKS_RLDC}
Sengupta T. K., Singh H., Bhaumik S., \& Chowdhury R. R., Diffusion in inhomogeneous flows: Unique equilibrium state in an internal flow, {\em Computers \& Fluids}, \textbf{88}, 440-451 (2013)

\bibitem{4_Buaria}
Doering C. R., The 3D Navier-Stokes Problem, {\em Ann. Rev. Fluid Mech.}, \textbf{41}, 109-128 (2009)

\bibitem{TKS_PRE12}
Sengupta T. K., Bhaumik S., \& Bhumkar Y. G., Direct numerical simulation of two-dimensional wall-bounded turbulent flows from receptivity stage, {\em Phys. Rev. E}, \textbf{85}(2), 026308 (2012)

\bibitem{TKS_IUTAM}
Sengupta T. K., A critical assessment of Simluations for Transitional and Turbulent Flows, {in:\em IUTAM Symp. Proc. Advances In Computation, Modeling And Control Of Transitional And Turbulent Flows}, 491-532 (2016)

\bibitem{TKS_vonNeumann07}
Sengupta T. K., Dipankar A., \& Sagaut P., Error dynamics: Beyond von Neumann analysis, {\em J. Comp. Phys.}, \textbf{226}(2), 1211-1218 (2007)

\bibitem{B74}
Sengupta T. K., \& Bhole A., Error dynamics of diffusion equation: Effects of numerical diffusion and dispersive diffusion, {\em J. Comp. Phys.}, \textbf{266}, 240-251 (2014)

\bibitem{B94}
Suman V. K., Sengupta T. K., Jyothi Durga Prasad C., Surya Mohan K., \& Sanwalia D., Spectral analysis of finite difference schemes for convection diffusion equation, {\em Computers \& Fluids}, \textbf{150}, 95-114 (2017)

\bibitem{David}
David Cl., Sagaut P., \& Sengupta T., A linear dispersive mechanism for numerical error growth: spurious caustics, {\em Euro. J. Mech. - B/Fluids}, \textbf{28}(1), 146-151 (2009)

\bibitem{Tan2021}
Tan R., Ooi A., \& Sandberg R. D., Two Dimensional Analysis of Hybrid Spectral/Finite Difference Schemes for Linearized Compressible Navier--Stokes Equations, {\em J. Sci. Comp.}, \textbf{87}(2), 42 (2021)

\bibitem{B118}
Suman V. K., Sengupta T. K., \& Mathur J. S., Effects of numerical anti-diffusion in closed unsteady flows governed by two-dimensional Navier-Stokes equation, {\em Computers \& Fluids}, \textbf{201}, 104479 (2020)

\bibitem{HACM}
Sengupta T. K., High Accuracy Computing Methods: Fluid Flows and Wave Phenomena, Cambridge Univ. Press, USA (2013)

\bibitem{20_Buaria}
Moin P., \& Mahesh K., Direct Numerical Simulation: A Tool in Turbulence Research, {\em Ann. Rev. Fluid Mech.}, \textbf{30}, 539-578 (1998)

\bibitem{5_Buaria}
Beale J. T., Kato T., \& Majda A., Remarks on the breakdown of smooth solutions for the 3-D Euler equations, {\em Commun. Math. Phys.}, \textbf{94}, 61-66 (1984)

\bibitem{vonNeumann_Richtmeyer}
von Neumann J., \& Richtmyer R. D., On the numerical solution of partial differential equations of parabolic type, {\em Los Alamos Rept., Series A}, \textbf{LA-657}, 1-17 (1947)

\end{thebibliography}
\end{document}